\documentclass{amsart}
\usepackage{amssymb}
\usepackage[all]{xy}
\CompileMatrices
\newtheorem*{introthm}{Theorem}
\newtheorem{thm}{Theorem}[section]
\newtheorem{prop}[thm]{Proposition}
\def\quot{/\!\!/}
\def\mal{\! \cdot \!}
\def\stern{\! * \!}
\def\KK{{\mathbb K}}
\def\CC{{\mathbb C}}
\def\AA{{\mathbb A}}
\def\ZZ{{\mathbb Z}}
\def\QQ{{\mathbb Q}}
\def\Aut{\operatorname{Aut}}

\def\cone{{\rm cone}}

\begin{document}
\title[Demushkin's Theorem in codimension one]
      {Demushkin's Theorem in codimension one}
\author[F. Berchtold]{Florian Berchtold} 
\address{Fachbereich Mathematik und Statistik, Universit\"at Konstanz,
         78457 Konstanz, Germany}
\email{Florian.Berchtold@uni-konstanz.de}
\author[J.~Hausen]{J\"urgen Hausen} 
\address{Fachbereich Mathematik und Statistik, Universit\"at Konstanz,
         78457 Konstanz, Germany}
\email{Juergen.Hausen@uni-konstanz.de}
\begin{abstract}
Demushkin's Theorem says that any two toric structures on an affine
variety $X$ are conjugate in the automorphism group of $X$. We provide
the following extension: Let an $(n-1)$-dimensional torus $T$ act
effectively on an $n$-dimensional affine toric variety $X$. Then $T$
is conjugate in the automorphism group of $X$ to a subtorus of the big
torus of $X$.
\end{abstract}

\maketitle

\section*{Introduction}

This paper deals with automorphism groups of toric varieties $X$
over an algebraically closed field $\KK$ of characteristic zero. 
We consider the following problem: Let $T \times X \to X$ be an 
effective regular torus action. When is this action conjugate in
$\Aut(X)$ to the action of a subtorus of the big torus $T_{X}
\subset X$? Some classical results are: 
\begin{itemize}
\item For complete $X$, the answer is always positive, because then
  $\Aut(X)$ is an affine algebraic group with maximal torus
  $T_{X}$, compare~\cite{Dem} and~\cite{Co}.
\item For $X = \KK^{m}$ and $\dim(T) \ge m-1$, positive answer is due
  to Bia\l ynicki-Birula, see~\cite{BB0} and~\cite{BB1}.
\item For $X$ affine and $\dim(T) = \dim(X)$, positive answer
  is due to Demushkin~\cite{De} and Gubeladze~\cite{Gu}.
\end{itemize}

We focus here on the case $\dim(T) = \dim(X)-1$. As in~\cite{De}
and~\cite{Gu}, we shall assume that $X$ has no torus factors. 
We do not insist on $X$ being affine; we just require that $X$ has no
``small holes'' in the sense that there is no open toric embedding $X
\to X'$ with $X' \setminus X$ nonempty of codimension at least
two. Under these assumptions we prove, see Theorem~\ref{result}:

\begin{introthm}
Let $T \times X \to X$ be an effective regular action of an 
algebraic torus $T$ of dimension $\dim(X)-1$. Then $T$ is conjugate in
$\Aut(X)$ to a subtorus of the big torus $T_{X} \subset X$.
\end{introthm}

In the case of tori $T$ of dimension strictly less than $\dim(X)-1$
the ``toric linearization problem'' stated at the beginning is wide
open, even for actions of an $(n-2)$-dimensional torus on the
affine space $\AA^{n}$. In the latter setting, there is a positive 
result for the case of a fixed point set of positive dimension,
see~\cite{KoRu1}, and a deep theorem saying that $\CC^{*}$-actions on
$\CC^{3}$ are linearizable,
see~\cite{Kaetal} and~\cite{KoRu2}.

Let us outline the main ideas of the proof of our theorem. 
In contrast to~\cite{De} and~\cite{Gu}, our approach is geometric. 
Since any two toric structures on $X$ are conjugate in the
automorphism group of $X$, see~\cite{Be}, it suffices to extend the
$T$-action to an almost homogeneous torus action on $X$. This is done
in three steps:

First lift the $T$-action (up to a finite homomorphism $T \to T$) 
to Cox's quotient presentation 
$\xymatrix{{\KK^{m}} \ar@{.>}[r] & X}$, see
Section~\ref{section1}. Next extend the lifted $T$-action to a toric
structure on $\KK^{m}$. This involves linearization of a certain
diagonalizable group action, see Section~\ref{section2}. Finally,
push down the new toric structure of $\KK^{m}$ to $X$. For this we
need that $X$ has no small holes, see Section~\ref{section3}.

\section{Lifting torus actions}
\label{section1}

We provide here a lifting result for torus actions on a toric variety
$X$ to the quotient presentation of $X$ introduced by Cox~\cite{Co}. 
First we recall the latter construction. For notation and the basic
facts on toric varieties, we refer to Fulton's book~\cite{Fu}. 

We shall assume that the toric variety $X$ is {\em nondegenerate},
that is $X$ admits no toric decomposition $X \cong Y \times \KK^{*}$. 
Note that this is equivalent to requiring that every invertible 
$f \in \mathcal{O}(X)$ is constant. 

Let $X$ arise from a fan $\Delta$ in a lattice $N$. Denote the rays of
$\Delta$ by $\varrho_{1}, \ldots, \varrho_{m}$. Let $Q \colon
\ZZ^{m} \to N$ be the map sending the canonical base vector $e_{i}$ to
the primitive generator of $\varrho_{i}$. For a maximal cone 
$\tau \in \Delta$, set
$$ \sigma(\tau) := \cone(e_{i}; \; \varrho_{i} \subset \tau). $$  

Then these cones $\sigma(\tau)$ are the maximal cones of a fan
$\Sigma$ consisting of faces of the positive orthant in $\QQ^{m}$. 
Moreover, $Q \colon \ZZ^{m} \to N$ is a map of the fans $\Sigma$ and
$\Delta$. The following properties of this construction are well
known:

\begin{prop}\label{cox}
Let $Z \subset \KK^{m}$ be the toric variety defined by
$\Sigma$, let $q \colon Z \to X$ be the toric morphism
corresponding to $Q \colon \ZZ^{m} \to N$, and let 
$H \subset T_{Z}$ be the kernel of the homomorphism 
$T_{Z} \to T_{X}$ of the big tori obtained by restricting 
$q \colon Z \to X$. 
\begin{enumerate}
\item The complement $\KK^{m} \setminus Z$ is of dimension
  at most $m-2$.
\item The map $q \colon Z \to X$ is a good quotient for the
  action of $H$ on $Z$. 
\item $X$ is smooth if and only if the group $H$ acts freely.
\end{enumerate}
\end{prop}

In general, the diagonalizable group $H \subset T_{Z}$ may be
disconnected. Hence we can at most expect liftings of a given action
$T \times X \to X$ in the sense that $q \colon Z \to X$ becomes 
$T$-equivariant up to a (finite) epimorphism $T \to T$. But such
liftings exist:
 
\begin{prop}\label{lifttocox}
Notation as in~\ref{cox}.
Let $T \times X \to X$ be an effective algebraic torus action. Then
there exist an effective regular action $T \times Z \to Z$ and an
epimorphism $\kappa \colon T \to T$ such that
\begin{enumerate}
\item $t \mal (h \mal z) = h \mal (t \mal z)$ holds for all $(t,h,z)
  \in T \times H \times Z$,
\item $q(t \mal z) = \kappa(t) \mal q(z)$ holds for all $(t,z) \in T
  \times X$.
\end{enumerate} 
\end{prop}

\proof First we reduce to the case that $X$ is smooth. So, suppose for
the moment that the assertion is proven in the smooth case. Then we can
lift the $T$-action over the set $U \subset X$ of smooth points. The
task then is to extend the lifted action from $U' := q^{-1}(U)$ to $Z$.

By Sumihiro's Theorem~\cite[Cor.~2]{Su}, $X$ is covered by $T$-invariant
affine open subsets $V \subset X$. The inverse images $V' :=
q^{-1}(V)$ are affine and $V' \setminus U'$ is of codimension at least
2 in $V'$. This allows to extend uniquely the lifted $T$-action from 
$V' \cap U'$ to $V'$ and hence from $U'$ to $Z$.

Therefore we may assume in the remainder of this proof that the toric
variety $X$ is smooth. As noted in Proposition~\ref{cox}, this means
that the group $H \subset T_{Z}$ acts freely on $Z$.
  
The most convenient way to lift the $T$-action is to split the
procedure into simple steps. For this, write $H$ as a direct product
of a torus $H_{0}$ with finite cyclic groups $H_{1}, \ldots, H_{k}$.
This gives rise to a decomposition of the quotient presentation 
$Z \to X$:
$$ 
\xymatrix{
Z \ar[r]^{/H_{k}} & 
Z_{k-1} \ar[r]^{/H_{k-1}} & 
& \cdots & 
{\ar[r]^{/H_{1}}} & 
Z_{0} \ar[r]^{/H_{0}} & 
X }
$$

This decomposition allows us to lift the $T$-action step by step with
respect to the geometric quotients by the free actions of the factors
$H_{i}$. We shall write again $Z$, $H$ and $X$ instead of
$Z_{i}$, $H_{i}$ and $Z_{i-1}$.

The action of $H$ on $Z$ defines a grading of the 
$\mathcal{O}_{X}$-algebra $\mathcal{A} := q_{*}(\mathcal{O}_{Z})$. 
Namely, denoting by $\Gamma$ the character group of $H$, we have
for every open $V \subset X$ the decomposition into
homogeneous functions: 
$$ 
\mathcal{O}(q^{-1}(V))
=
\mathcal{A}(V) 
= 
\bigoplus_{\chi \in \Gamma} \mathcal{A}_{\chi}(V). 
$$

Since $H$ acts freely on $Z$, all homogeneous components
$\mathcal{A}_{\chi}$ are locally free $\mathcal{O}_{X}$-modules of
rank one. We shall use this fact to make the $\mathcal{O}_{X}$-algebra 
$\mathcal{A}$ into a $T$-sheaf over the $T$-variety $X$. Then it is
canonical to extract the desired lifting from this $T$-sheaf structure. 

If the group $H$ is connected, then we can prescribe
$T$-linearizations on the $\mathcal{O}_{X}$-modules $\mathcal{A}_{i}$
corresponding to the members $\chi_{i}$ of some lattice basis of
$\Gamma$. Tensoring these linearizations gives the desired $T$-sheaf
structure on the $\mathcal{O}_{X}$-algebra $\mathcal{A}$, compare
also~\cite[Section~3]{Ha}.
 
Since $X$ is covered by $T$-invariant affine open subsets, we can
easily check that this $T$-sheaf structure of $\mathcal{A}$ arises
from a regular $T$-action on $Z$ that commutes with the action of $H$
and makes the quotient map $q \colon Z \to X$ even equivariant. This
settles the case of a connected $H$. 

Assume that $H$ is finite cyclic of order $d$. Let $\chi$
be a generator of $\Gamma$. Again, we choose a $T$-linearization
of $\mathcal{A}_{\chi}$. But now it may happen that the induced
$T$-linearization on $\mathcal{A}_{d \chi} = \mathcal{O}_{X}$ is not
the canonical one. However, since $\mathcal{O}^{*}(X) = \KK^{*}$
holds, these two linearizations only differ by a character $\xi$ of
$T$.

Let $\kappa \colon T \to T$ be an epimorphism such that 
$\xi \circ \kappa = \xi_{0}^{d}$ holds for some character $\xi_{0}$ of
$T$. Consider the action $t \stern x := \kappa(t) \mal x$ on $X$.
Then $\mathcal{A}_{\chi}$ is also linearized with respect to this
action by setting $t \stern f := \kappa(t) \mal f$. Twisting with
$\xi_{0}^{-1}$, we achieve that the induced linearization on 
$\mathcal{A}_{d \chi} = \mathcal{O}_{X}$ is the canonical one:
$$ 
(t \stern f)(x) 
= 
\xi_{0}^{-d}(t) \xi(\kappa(t)) f(t^{-1} \stern x)
=
f(t^{-1} \stern x).
$$

The rest is similar to the preceding step: The $T$-sheaf structure of
$\mathcal{A}$ defines a $T$-action $(t,z) \mapsto t \stern z$ on $Z$
commuting with the action of $H$ and making the quotient map $q \colon
Z \to X$ equivariant with respect to $(t,x) \mapsto t \stern x$. 
Dividing by the kernel of ineffectivity, we can make the action on $Z$
effective and obtain the desired lifting. \endproof

\section{Diagonalizable group actions}
\label{section2}

In this section we show that any effective regular action of an
$(m-1)$-dimensional diagonalizable group $G$ on $\KK^{m}$ can be
brought into diagonal form by means of an algebraic coordinate
change. The result extends a well known analogous statement on torus
actions due to Bia\l ynicki-Birula, see~\cite{BB1}. 

We would like to thank the referee for his valuable proposals in order
to make our first proof more transparent.  

\begin{prop}\label{diagdiag}
Let $G \times \KK^{m} \to \KK^{m}$ be an effective algebraic action of
an $(m-1)$-dimensional diagonalizable group $G$. Then there exist
$\alpha \in \Aut(\KK^{m})$ and characters  $\chi_{i} \colon G \to
\KK^{*}$ such that for every $g \in G$ and every $(z_{1}, \ldots,
z_{m}) \in \KK^{m}$ we have
$$ 
\alpha(g \mal \alpha^{-1}(z_{1}, \ldots, z_{m})) 
= 
(\chi_{1}(g)z_{1}, \ldots, \chi_{m}(g)z_{m}).
$$
\end{prop}

\proof If the quotient space $\KK^{m} \quot G$ is a point, then an
application of Luna's slice theorem shows that the action of $G$ is
linearizable; this works even more generally for reductive groups,
see~\cite[Proposition~5.1]{KrPo}. Since any linear $G$-action is 
diagonalizable, we obtain the assertion in the case of 
$\KK^{m} \quot G$ being a point.

So we are left with the case that $\KK^{m} \quot G$ is of dimension
one. We write $G = G_{0} \times G_{1}$ with an algebraic 
torus $G_{0}$ and a finite abelian group $G_{1}$. According to the
main result of~\cite{BB1}, we may assume that the action of $G_{0}$ is
already diagonal. In the sequel, we view $G_{0}$ as a subtorus of
the torus $(\KK^{*})^{m}$. 

We shall show that the group $G_{1}$ permutes the coordinate
hyperplanes $V(z_{i})$ of $\KK^{m}$. Indeed, this is all we need,
because then $G_{1}$ acts by linear automorphisms and hence the action
of $G$ is diagonalizable.  

Let $p \colon \KK^{m} \to \KK^{m} \quot G_{0}$ be the quotient 
map. This is a toric morphism. In particular, since the quotient space
$\KK^{m} \quot G_{0}$ is of dimension one, it is isomorphic to $\KK$. 
Thus we can write down the quotient map explicitly:
There are relatively prime nonnegative integers $a_{i}$ such that 
$$ 
p(z_{1}, \ldots, z_{m}) = z_{1}^{a_{1}} \ldots z_{m}^{a_{m}}.
$$

Let us renumber the coordinates such that in the above presentation
of the quotient map we have $a_{i} > 0$ for all $i \le k$ and $a_{i} =
0$ for all $i > k$ with a suitable integer $k \le m$. 

Consider an $i > k$, say $i = m$. Then the points $(1, \ldots, 1)$
and $(1,\ldots,1,0)$ lie in the regular fiber $p^{-1}(1)$ of the
quotient map $p \colon \KK^{m} \to \KK$. Thus $(1,\ldots,1,0)$ lies in
the closure of $G_{0}$ and hence it lies in the closure of some
onedimensional subtorus $T_{m} \subset G_{0}$.  
But this subtorus is necessarily of the form
$$T_{m} = \{(1,\ldots,1,t); \; t \in \KK^{*} \}. $$

Since the actions of $G_{1}$ and $T_{m}$ commute, the group
$G_{1}$ leaves the fixed point set of $T_{m}$ invariant. But the
latter is just the coordinate hyperplane $V(z_{m})$. Analogously, we
conclude that the remaining $V(z_{i})$, where $i>k$, are invariant
under the group $G_{1}$. 

We discuss now what happens to the coordinate hyperplanes
$V(z_{i})$, where $i \le k$. Note that these are precisely the
irreducible components of the fiber $p^{-1}(0)$ of the quotient
map. We shall distinguish the cases $k>1$ and $k=1$. 

For $k>1$, the fiber $p^{-1}(0)$ is the unique reducible fiber
of the quotient map. On the other hand the action of $G_{1}$ commutes
with the action of $G_{0}$ and hence $G_{1}$ permutes the fibres of
$p$. Thus $G_{1}$ has to leave $p^{-1}(0)$ invariant. Hence $G_{1}$
permutes the coordinate hyperplanes $V(z_{1}), \ldots, V(z_{k})$
provided $k > 1$. 

Finally, we treat the case $k=1$. Then, as seen before, $G_{0}$ 
contains the one dimensional tori $T_{2}, \ldots, T_{m}$. In other
words, we have 
$$ G_{0} = \{(1,t_{2}, \ldots, t_{m}); \; t_{i} \in \KK^{*}\}. $$ 

So the fixed point set of $G_{0}$ is just the $z_{1}$-axis. Hence 
$G_{1}$ leaves the $z_{1}$-axis invariant. Now, $G_{1}$ acts on 
the $z_{1}$-axis with a fixed point, say $b$. By conjugating the
$G$-action with the translation by $b$, we achieve that $b=0$
holds. But then $G_{1}$ leaves $V(z_{1}) = p^{-1}(0)$
invariant. \endproof

\section{Proof of the main result}
\label{section3}

We say that a toric variety $X$ {\em has no small holes}, if it does not
admit an open toric embedding $X \subset X'$ such that $X' \setminus
X$ is nonempty of codimension at least $2$ in $X'$. Examples are
the toric varieties arising from a fan with convex support. This
comprises in particular the affine ones. 

\begin{thm}\label{result}
Let $X$ be a nondegenerate toric variety without small holes, and let
$T \times X \to X$ be an effective regular action of an algebraic
torus $T$ of dimension $\dim(X)-1$. Then $T$ is conjugate in $\Aut(X)$
to a subtorus of the big torus $T_{X} \subset X$.  
\end{thm}

\proof According to~\cite[Theorem~4.1]{Be}, any two toric structures
of $X$ are conjugate in $\Aut(X)$. Consequently, it suffices to
show that the action of $T$ on $X$ extends to an effective regular
action of a torus of dimension $\dim(X)$ on $X$. 

Consider Cox's construction $q \colon Z \to X$ and its
kernel $H := \ker(q)$ as defined in~\ref{cox}. Choose a lifting of the
$T$-action to $Z$ as provided by Proposition~\ref{lifttocox}. This
gives us an action of the $(m-1)$-dimensional diagonalizable group $G
:= T \times H$ on the open set $Z \subset \KK^{m}$. 

Since the complement $\KK^{m} \setminus Z$ is of dimension at most
$m-2$, the action of $G$ extends regularly to $\KK^{m}$. Let $G_{0}$
be the (finite) kernel of ineffectivity. Applying
Proposition~\ref{diagdiag} to the action of $G/G_{0}$, we can extend
the $G$-action to an almost homogeneous action of a torus $S$ on
$\KK^{m}$. 

We show that $Z$ is invariant with respect to the action of $S$. 
According to~\cite[Corollary~2.3]{Sw}, we obtain this if we can 
prove that the set $Z$ is $H$-maximal in the following sense: If $Z'
\subset \KK^{m}$ is an $H$-invariant open subset admitting a good
quotient $q' \colon Z' \to X'$ by the action of $H$ such that $Z$ is a
$q'$-saturated open subset of $Z'$, then we already have $Z' = Z$. 

To verify $H$-maximality of $Z$, consider $Z' \subset \KK^{m}$ and $q'
\colon Z' \to X'$ as above. We may assume that $Z'$ is
$H$-maximal in $\KK^{m}$. Applying~\cite[Corollary~2.3]{Sw} to the
actions of $H$ and the standard torus $(\KK^{*})^{m}$ on $\KK^{m}$, we
obtain that $Z'$ is invariant with respect to the action of
$(\KK^{*})^{m}$. Hence we obtain a commutative diagram of toric
morphisms:
$$ 
\xymatrix{
Z \ar[r] \ar[d]^{q}_{\quot H} 
& 
Z' \ar[d]_{q'}^{\quot H} \\
X \ar[r] 
& 
X'}
$$

By the choice of $Z'$, the horizontal arrows are open
toric embeddings. Moreover, the complement $X' \setminus X$ is of
codimension at least two in $X'$, because its inverse image
$Z'\setminus Z$ under $q'$ is a subset of the small set 
$\KK^{m} \setminus Z$ and $q'$ is surjective. By the assumption on
$X$, we obtain $X' = X$. This verifies $H$-maximality of $Z$. Hence
our claim is proved. 

The rest is easy: The torus $S/H$ acts with a dense orbit on
$X$. Dividing $S/H$ by the kernel of ineffectivity of this action, 
we obtain the desired extension of the action of $T$ on $X$. \endproof

\end{document}